\input amstex
\documentstyle{amsppt}

\topmatter

\title
\nofrills
          On the Fay identity for KdV tau  functions and 
           the identity for the Wronskian of squared solutions 
                      of Sturm-Liouville equation
\endtitle

\author 
          Yordan P. Mishev 
\endauthor

\rightheadtext{Fay identity and squared solutions ...}

\affil
          Research Institute for Mathematical Sciences, Kyoto University, \\
                             Sakyo-ku 606, Kyoto, JAPAN \\
               e-mail  address : mishevyp\@kurims.kyoto-u.ac.jp \\
           On leave from: Forestry University, Sofia, BULGARIA.
\endaffil

\date 
          March, 1998
\enddate

\abstract
We show that the well known  identity for the Wronskian 
of squared solutions of a Sturm-Liouville equation follows from the  
Fay identity.
We also study some odd-order ( ($2^n -1$)-order , $ n = 2 , 3 ,  \dots$ ) 
identities which are specific for tau functions, related to the KdV hierarchy. 
\endabstract

\endtopmatter

\document 

\subhead
I. Introduction
\endsubhead

We began this research as a study of the expression of the
Wronskian of squared solutions of Sturm-Liouville equation 
by KdV tau functions and Fay identity.    
Now, when the desired result is obtained (cf. Theorem 1.1), we realize that 
this is a story of the common origine of the following three relations for 
the functions: $\text{x}, \, \sin (x),\, \theta_{11} (x)$ respectivelly 
( $  x, z_1, z_2 \in \Bbb C $ ):
$$
\align
&(z_2-z_1)\biggl[ (x+z_1+z_2)(x-z_1)(x-z_2)-
                  (x-z_1-z_2)(x+z_1)(x+z_2)\biggr]  =  \\ 
&(z_1+z_2)\biggl[ (x+z_1-z_2)(x-z_1)(x+z_2)-(x-z_1+z_2)(x+z_1)(x-z_2)\biggr] , 
\endalign
$$
$$
\align
\sin(z_2-z_1)\biggl[&\sin(x+z_1+z_2)\sin(x-z_1)\sin(x-z_2)- \\
                    &\sin(x-z_1-z_2)\sin(x+z_1)\sin(x+z_2)\biggr] =  \\ 
\sin(z_1+z_2)\biggl[&\sin(x+z_1-z_2)\sin(x-z_1)\sin(x+z_2)- \\ 
                    &\sin(x-z_1+z_2)\sin(x+z_1)\sin(x-z_2)\biggr] ,
\endalign
$$
$$
\align
\theta_{11}(z_2-z_1)\biggl[&\theta_{11}(x+z_1+z_2)
                            \theta_{11}(x-z_1)\theta_{11}(x-z_2)- \\
&\theta_{11}(x-z_1-z_2)\theta_{11}(x+z_1)\theta_{11}(x+z_2)\biggr]=  \\ 
\theta_{11}(z_1+z_2)\biggl[&\theta_{11}(x+z_1-z_2)
                            \theta_{11}(x-z_1)\theta_{11}(x+z_2)- \\
   &\theta_{11}(x-z_1+z_2)\theta_{11}(x+z_1)\theta_{11}(x-z_2)\biggr] 
\endalign
$$
(we use the notations for theta functions from Ref. 1.

In the present paper we will prove an identity for general KdV tau 
functions (it will be a trird-order identity for tau). 
The mentioned three types of functions are roughly speaking three 
types of KdV tau functions for stationary ($t$-independent) solutions 
of the KdV equation:
$
u_t = 6 u u_x + u_{xxx}.
$ 
So, the polynomial relation follows from this cubic identity for KdV taus, but
there are some specific problems to translate the identity
for the general KdV taus to the cases of trigonometric functions and
elliptic theta functions. We will postpone the solution of these problems  
to some next publication.  

This way, the status of the three relations is quite different: the first one 
is easy to prove directly (and it also follows from the cubic identity 
for KdV taus); the second one is not difficult to prove directly, using 
the well known trigonometric identities; and the third one is still 
conjectural (we could not derive it from the Riemann relations, 
we could only check it numerically, using the system Mathematica 3.0
at RIMS, Kyoto University). 

Let $\tau(t)$ , $t \equiv ( t_1, t_2, t_3,\dots ) \in \Bbb C ^{\infty}$ ,  
$t_1 \equiv x $ is an arbitrary  tau function, related to 
the Kadomtzev-Petviashvili (KP) hierarchy Ref. 2.
Let us denote ($z \in \Bbb C$): 
$$
\align
[z] &:= (z, z^2/2, z^3/3, \dots ) \in \Bbb C^{\infty}, \\
\tau (t+[z]) &:= \tau (t_1 + z, t_2 + z^2/2, t_3 + z^3/3, \dots ) .
\endalign
$$ 

The following identity ($z_0 , z_1 , z_2 , z_3 \in \Bbb C$) :
$$
\align
(z_0-z_1)(z_2-z_3) \tau (t+[z_0]+[z_1]) \tau (t+[z_2]+[z_3]) &+ \\
(z_0-z_2)(z_3-z_1) \tau (t+[z_0]+[z_2]) \tau (t+[z_3]+[z_1]) &+ \tag 1.1 \\
(z_0-z_3)(z_1-z_2) \tau (t+[z_0]+[z_3]) \tau (t+[z_1]+[z_2]) &= 0
\endalign
$$
is called {\it Fay identity} Ref. 3 for the KP tau function $\tau$.
It was firstly obtained Ref.4 for theta functions related to  
Jacobians. In genus $g =1$ case its form is:
$$
\align
\theta_{11} (z_0-z_1) \theta_{11} (z_2-z_3) 
                      \theta_{11} (t+z_0+z_1) \theta_{11} (t+z_2+z_3) &+ \\
\theta_{11} (z_0-z_2) \theta_{11} (z_3-z_1) 
                      \theta_{11} (t+z_0+z_2) \theta_{11} (t+z_3+z_1) &+  \\
\theta_{11} (z_0-z_3) \theta_{11} (z_1-z_2) 
                      \theta_{11} (t+z_0+z_3) \theta_{11} (t+z_1+z_2) &= 0
\endalign
$$ 
Afterwards it was used Ref. 1 in geometric tratement of 
soliton equations. 
Leter it was generalized for tau functions Ref. 3.

Fay identity is fulfill also for tau functions related to $n$ - th 
( $n$ = 2, 3, 4, ... ) Gel'fand-Dickey reduction of KP hierarchy. In 
the present paper we will consider only the $n = 2$ reduction, i.e.
the KdV hierarchy. Such tau functions we will call 
KdV tau functions. They can be characterized by the conditions 
($\partial_{t_{2k}} \equiv  {\partial}/{\partial t_{2k}} $):
$$
   \partial_{t_{2k}} \, \tau(t) = 0 \, ,
   \qquad  k = 1, 2, 3, \dots
$$ 
which imply for every $z \in \Bbb C$ : 
$$
   \tau (t-[z]) = \tau (t+[-z]) .  \tag 1.2
$$
 
There are two main goals in the present article.  
The first aim is 
to show that the famous identity for the Wronskian 
$\bigl(  W (f, g) := fg'-f'g$ , 
$' \equiv \partial_x \equiv \frac {\partial} {\partial x} \bigr)$
of squared solutions of Sturm-Liouville equation Ref. 5 
follows from the  Fay identity for $KdV$ tau functions. 
The second aim is to obtain 
some specific relations for the KdV tau functions. 

We came to these results when studying the problem  
to find a dictionary between the tau functions and some  
formulas related to squared solutions of Sturm-Liouville equation 
(especially the mentioned identity for the Wronskian of 
squared solutions- an important ingredient of this area 
Ref. 5, Ref. 6. 
Such dictionary will be useful in examining some features of Miura 
transformations. 
It is well known that squared solutions span 
the kernels of the Frechet derivatives
$\,  \frak M '_{\pm} (v) = 2 v \pm \partial_x \,$
of Miura transformations 
$\, u_{\pm} = \frak M_{\pm} v := v^2 \pm  v_x \,$ , 
($\, v_x \equiv \partial_x v \, $ ), where
$\, v \,$ is a solution of the mKdV equation and $\, u_{\pm} \,$  
are solutions of the KdV equation. 
It is also well known Ref. 7 the interpretation of 
Miura transformations as projections from flag to corresponding subspace 
(in Sato Grassmannian), which is intimately connected to tau functions. 
Some parts of the dictionary were known: e.g. a formula, which expresses 
the squared solutions by means of $\tau (t)$ and vertex operators 
Ref. 2 (in this paper the so-called 
$\Lambda$-operators Ref. 5 are also mentioned). 
There were no relations to the tau functions of the identity 
for the Wronskian of squared solutions, 
but it was well known an expression of the Wronskian of two
solutions by means of $\tau (t)$ Ref. 2, Ref. 8. 
Because of the fact that in the proof of the latest formula, the Fay identity 
was used, we expected that the same identity will be useful in the 
``paraphrase'' of the Wronskian of equared solutions.

We need such dictionary, because we observed some similarities between 
Matsuo and Cherednik transformations (Ref. 9):  
$$ 
\text {Knizhniv-Zamolodchikov equation}       \to
\text {quantum Calogero-Sutherland  system}
$$
on the one hand and Miura transformation  
on the other hand. Our oppinion is that such similarities will be 
easier explained on the language of tau functions, flag and Grassmann 
manifolds, etc.. So, the ``paraphrase'' of the relations for the 
Wronskian of equared solutions of Sturm-Liouville equation is only 
the first step in this direction.
We also expect that the presented connections between squared solutions 
and tau functions will be useful in another areas of the subject 
(cf. Ref.10, Ref. 11). 

In order to explain the main results of the present article, let us 
remind some notations Ref. 2, Ref. 8. Let $\psi (x, z)$ and 
$\psi ^{*}(x, z)$ are two linearly independent solutions (cf. Section II)
of the Sturm-Liouville equation :
$$
\bigl( \partial^{2}_{x} + u(x) \bigr) \, \psi(x,z) = z^2 \psi(x,z) . \tag 1.3 
$$
Then the following relations : 
$$
W \bigl( \psi (x, z_1) \psi^{*} (x, z_1) , 
         \psi (x, z_2) \psi^{*} (x, z_2) \bigr) = 
$$
$$
\align
  - (z_1^2-z_2^2)^{-1} \, \partial_{x} &\biggl[ W \bigl( 
  \psi (x, z_1) , \psi (x, z_2) \bigr)  W \bigl( \psi^{*} (x, z_1), 
       \psi^{*} (x, z_2) \bigr) \biggr] =  \\  
   - (z_1^2-z_2^2)^{-1} \, \partial_{x} &\biggl[ W \bigl( 
  \psi (x, z_1) , \psi  ^{*}(x, z_2) \bigr) W \bigl( \psi^{*} (x, z_1), 
       \psi (x, z_2) \bigr) \biggr] ,  \tag 1.4  
\endalign 
$$
($z_1, z_2 \in \Bbb C$) we will call {\it Faddeev-Tahtajan identity}. 
This relation has a long history. It was used in the theory of
inverse spectral problems for the Sturm-Liouville operators. Afterwards 
Faddeev-Tahtajan identity 
played an important role in the first years of Soliton Theory. In 
Ref. 5 the origin of the identity is interpreted in terms of  
classical $r$-matrixes. Here we will explain the origin of 
Faddeev-Tahtajan identity using the language of tau functions.
 
The first main result in this paper is given in the following 
\proclaim {Theorem 1.1}
The Faddeev-Tahtajan identity (1.4) follows from the Fay identity 
(1.1) for KdV tau functions.
\endproclaim
 
The second main result in the present article is given in the following 
\proclaim{Theorem 1.2} 
Let $\tau (t)$ , $t \in \Bbb C ^{\infty}$ is an arbitrary KdV tau function.
Then :

(i) for every $ z_1, z_2 \in \Bbb C $ : 
$$
\multline
(z_2-z_1)\biggl[ \tau (t+[z_1]+[z_2]) \tau (t-[z_1]) \tau (t-[z_2]) -      
             \tau (t-[z_1]-[z_2]) \tau (t+[z_1]) \tau (t+[z_2]) \biggr]=  \\ 
\quad
  (z_2+z_1)\biggl[  \tau (t+[z_1]-[z_2]) \tau (t-[z_1]) \tau (t+[z_2]) -     
             \tau (t-[z_1]+[z_2]) \tau (t+[z_1]) \tau (t-[z_2]) \biggr]
\, ; 
\endmultline
$$

(ii) for every $z \in \Bbb C$ :
$$\multline
\tau (t+2[z]) \tau^2 (t-[z]) - \tau (t-2[z]) \tau^2 (t+[z]) = \\ 
\quad
2 \sum_{k=0}^{\infty} z^{2k+1} \biggl[ \tau (t-[z]) \,  W _{2k+1} 
                          \bigl( \tau (t), \tau (t+[z]) \bigr) +
                                       \tau (t+[z])  \, W _{2k+1} 
                          \bigl( \tau (t), \tau (t-[z]) \bigr) \biggr] ,
\endmultline 
$$

where we denote :
$
  W _{2k+1} \bigl( f, g \bigr) := f (\partial_{t_{2k+1}} \, g) 
  - (\partial_{t_{2k+1}} \, f) g , \, k = 0, 1, 2, \dots  \, .
$ 
\endproclaim

{\bf Remark 1.3} 
Let us mention that the identities from Theorem 1.2 are cubic in $\, \tau$
relations (in contrary to the  Fay identity, which is quadratic 
in $\, \tau$ relation) and they are specific only for the KdV tau functions.

{\bf Remark 1.4}
The proof of Theorem 1.2 is based {\it only} on the following three
facts: 

(i) The Fay identity (1.1) (which is common for all tau functions),

(ii) The relation (1.2) (which is specific only for KdV tau functions),

(iii) The obvious identity for Wronskians :
$$
\align
 W \bigl( f_1f_2, g_1g_2 \bigr) \, = \,  &f_1 g_1 W \bigl( f_2, g_2 \bigr) +
                                    f_2g_2  W \bigl( f_1, g_1 \bigr) \, = \\
                                  &f_1g_2  W \bigl( f_2, g_1 \bigr) + 
                                   f_2g_1  W \bigl( f_1, g_2 \bigr) . \tag 1.5 
\endalign
$$
 
The paper consists of four sections. In Section II we give some preliminary
results. The proofs of the Theorem 1.1 and Theorem 1.2 are given in 
Section III.
In Section IV we give some examples and comments of the main statements.
A preliminary (and from different viewpoint) version of some of the results 
is presented in Ref. 12 and Ref. 13.

\subhead
II. Preliminary results
\endsubhead
  
Firstly, let us mention some obvious relations for Wronskians.
\proclaim {Lemma 2.1}
$$
  W \bigl( e^{z_1 x} f \, , \, e^{z_2 x} g \bigr)  =  
  e^{(z_1+z_2)x}\bigl[ \, W (f , g)  - (z_1-z_2) \, fg \, \bigr] , \tag i
$$
$$
  W \bigl( \frac {f_1}{g}, \frac {f_2}{g} \bigr)  = 
  \frac {W ( f_1 , f_2 )}{g^2} , \tag ii
$$
$$
  \partial _{x} \bigl( \frac {f_1 \, f_2}{g^2} \bigr) = 
  - \, \frac { f_1 W(f_2,g) + f_2 W(f_1,g) }{g^3} .  \tag iii
$$
\endproclaim 
Instead of Fay identity (1.1) we will use the  differential 
Fay identity Ref. 8 \linebreak  ($ z_1, z_2 \, \in \Bbb C$):
$$
\align
W \bigl( \tau (t+[z_1]) , \tau (t+[z_2]) \bigr) = (z_2^{-1}-z_1^{-1})
\bigl[ &\tau (t+[z_1])  \tau (t+[z_2]) - \\
       &\tau(t) \tau (t+[z_1]+[z_2]) \bigr] . \tag 2.1
\endalign
$$
Shifting the argument $\, t \,$ respectively to ($\, t-[z_1]-[z_2] \,$) ,  
($\, t-[z_2] \,$) and ($\, t-[z_1] \,$) we could obtain 
expressions respectively for the following Wronskians:
$$ W \bigl( \tau (t-[z_1] , \tau (t-[z_2]) \bigr) \, ,\quad 
   W \bigl( \tau (t+[z_1]-[z_2]) , \tau (t) \bigr) \, , \quad
   W \bigl( \tau (t-[z_1]+[z_2]) , \tau (t) \bigr) \, .
$$
But, shifting $\, t \,$ we cannot obtain an expression e.g. for the 
Wronskian : 
$$ W \bigl( \tau (t+[z_1]) , \tau (t-[z_2]) \bigr) . $$ 
This is possible for KdV tau functions. Using (2.1)
and (1.2), it is easy to see that: 
$$
\align
 W \bigl( \tau (t+[z_1])\, , \, \tau (t-[z_2]) \bigr) \, = \,
                                                 -(z_2^{-1}+z_1^{-1})
   \bigl[ &\tau (t+[z_1]) \tau (t-[z_2]) - \\ 
          &\tau (t) \tau (t+[z_1]-[z_2]) \bigr] .
\endalign
$$
This way we obtain the following expressions for the Wronskians of 
KdV tau functions.

\proclaim{Lemma 2.2}

Let $\, \tau (t) \, $ is an arbitrary KdV tau function. Then we have :
$$
\align
 W \bigl( \tau (t+[z_1])\, , \, \tau (t+[z_2]) \bigr) \, = \, 
                                                  (z_2^{-1}-z_1^{-1})
   \bigl[ &\tau (t+[z_1]) \tau (t+[z_2]) - \tag i \\
          &\tau(t) \tau (t+[z_1]+[z_2]) \bigr] , \\
 W \bigl( \tau (t-[z_1])\, , \, \tau (t-[z_2]) \bigr) \, = \,
                                                 -(z_2^{-1}-z_1^{-1})
   \bigl[ &\tau (t-[z_1]) \tau (t-[z_2]) - \\ 
          &\tau (t) \tau (t-[z_1]-[z_2]) \bigr] , \\
 W \bigl( \tau (t-[z_1])\, , \, \tau (t+[z_2]) \bigr) \, = \,
                                                  (z_2^{-1}+z_1^{-1}) 
   \bigl[ &\tau (t-[z_1]) \tau (t+[z_2]) - \\
          &\tau (t) \tau (t-[z_1]+[z_2]) \bigr] ,  \\      
 W \bigl( \tau (t+[z_1])\, , \, \tau (t-[z_2]) \bigr) \, = \,
                                                - (z_2^{-1}+z_1^{-1}) 
   \bigl[ &\tau (t+[z_1]) \tau (t-[z_2]) - \\
          &\tau (t) \tau (t+[z_1]-[z_2]) \bigr] ; \\ 
W \bigl( \tau (t+[z_1]-[z_2])\, , \, \tau (t) \bigr) \, = \,  
                                                  (z_2^{-1}-z_1^{-1})
   \bigl[ &\tau (t+[z_1]-[z_2]) \tau (t) - \tag ii \\
          &\tau (t+[z_1]) \tau (t-[z_2]) \bigr] \\
W \bigl( \tau (t-[z_1]+[z_2])\, , \, \tau (t) \bigr) \, = \, 
                                                - (z_2^{-1}-z_1^{-1}) 
   \bigl[ &\tau (t-[z_1]+[z_2]) \tau (t) - \\
          &\tau (t-[z_1]) \tau (t+[z_2]) \bigr] , \\     
W \bigl( \tau (t-[z_1]-[z_2])\, , \, \tau (t) \bigr) \, = \, 
                                                 (z_2^{-1}+z_1^{-1}) 
   \bigl[ &\tau (t-[z_1]-[z_2]) \tau (t) -  \\
          &\tau (t-[z_1]) \tau (t-[z_2]) \bigr] , \\     
W \bigl( \tau (t+[z_1]+[z_2])\, , \, \tau (t) \bigr) \, = \, 
                                                - (z_2^{-1}+z_1^{-1}) 
   \bigl[ &\tau (t+[z_1]+[z_2]) \tau (t) -  \\
          &\tau (t+[z_1]) \tau (t+[z_2]) \bigr].    
\endalign
$$
\endproclaim 

Let us define the {\it wave functions} $\psi (t,z)$ , $\psi ^{*} (t,z)$ 
Ref. 8 by expressions : 
$$
\align
\psi (t,z) &:= \exp(\sum_{k=1}^{\infty}{t_k z^k}) 
                                      \frac {\tau (t-[z^{-1}])}{\tau (t)} , \\ 
\psi ^{*} (t,z) &:= \exp(-\sum_{k=1}^{\infty}{t_k z^k}) 
                                      \frac {\tau (t+[z^{-1}])}{\tau (t)}.
\endalign
$$
For an arbitrary KdV tau function $\, \tau (t) \,$ ,
denoting $\, u(t) := 2 \, \partial_{x}^{2} \ln \tau (t) $ , it is 
well known Ref. 8 that the wave functions $\psi (t,z)$ and 
$\psi ^{*}(t,z)$ 
satisfy the Sturm-Liouville equation (1.3) 
($t_1 \equiv x$ and $t_3, t_5 ,\dots $ are parameters).
Using the relations of Lemma 2.2 we can explain the Wronskians
of the wave functions $\psi (t,z)$ and $\psi ^{*} (t,z)$ in terms of the
tau-function $\tau (t)$ .

\proclaim{Lemma 2.3}

Let $\, \tau (t) \, $ is an arbitrary KdV tau function and 
$\, \psi (t,z)$ , $\, \psi ^{*} (t,z) \, $ are the corresponding wave
functions. Then we have ($z_1 , z_2 \in \Bbb C$) : 

$$
\align
W \bigl( \psi (t,z_1) , \psi (t,z_2) \bigr)  =  (z_1-z_2) 
 \exp\bigl( \sum_{k=0}^{\infty} &t_{2k+1}(z_1^{2k+1}+z_2^{2k+1}) \bigr) 
                                                              \tag i \\
 &\frac {\tau (t-[z_1^{-1}]-[z_2^{-1}])}{\tau (t)} \, , \\ 
W \bigl( \psi ^{*} (t,z_1), \psi ^{*} (t,z_2) \bigr)  =  -(z_1-z_2)    
 \exp\bigl( -\sum_{k=0}^{\infty} &t_{2k+1}(z_1^{2k+1}+z_2^{2k+1})\bigr) 
                                                               \tag ii \\
 &\frac {\tau (t+[z_1^{-1}]+[z_2^{-1}])}{\tau (t)} \, , \\
W \bigl( \psi (t,z_1) , \psi ^{*} (t,z_2) \bigr)  =  (z_1+z_2) 
 \exp\bigl( \sum_{k=0}^{\infty} &t_{2k+1}(z_1^{2k+1}-z_2^{2k+1})\bigr)  
                                                                \tag iii \\
 &\frac {\tau (t-[z_1^{-1}]+[z_2^{-1}])}{\tau (t)} \, , \\
W \bigl( \psi ^{*} (t,z_1) , \psi (t,z_2) \bigr)  =  -(z_1+z_2)  
 \exp\bigl( -\sum_{k=0}^{\infty} &t_{2k+1}(z_1^{2k+1}-z_2^{2k+1})\bigr)  
                                                                 \tag iv \\
 &\frac {\tau (t+[z_1^{-1}]-[z_2^{-1}])}{\tau (t)} \, . 
\endalign  
$$
\endproclaim

{\it Proof :}

Let us denote the functions : 
$$
\varphi (t,z) := {e}^{zx} \frac {\tau (t-[z^{-1}])}{\tau (t)} \, , 
\qquad \qquad
\varphi ^{*} (t,z) := {e}^{zx} \frac {\tau (t+[z^{-1}])}{\tau (t)} \, .  
$$
Then we have :
$$
\psi (t,z) = \exp \bigl( \sum_{k=1}^{\infty} t_{2k+1} z^{2k+1} \bigr) 
                                           \varphi (t,z) , \qquad
\psi ^{*} (t,z) = \exp \bigl( - \sum_{k=1}^{\infty} t_{2k+1} z^{2k+1} \bigr) 
                                           \varphi ^{*} (t,z)  ,
$$
and consequently we have :
$$
W \bigl( \psi (t,z_1) , \psi (t,z_2) \bigr) = 
\exp \bigl( \sum_{k=1}^{\infty} t_{2k+1} (z_1^{2k+1} + z_2^{2k+1}) \bigr) \, 
W \bigl( \varphi (t,z_1) , \varphi (t,z_2) \bigr),\qquad  \text{etc.}
$$
Using the relations of Lemma 2.1 and Lemma 2.2 we obtain :
$$
W \bigl( \varphi (t,z_1) , \varphi (t,z_2) \bigr) = 
$$ 
$$ 
W \biggl( {e}^{z_1 x} \frac {\tau (t-[z_1^{-1}])}{\tau (t)} \, , \,  
          {e}^{z_2 x} \frac {\tau (t-[z_2^{-1}])}{\tau (t)}  \biggr) \, = 
$$
$$
{e}^{(z_1+z_2) x} \biggl[ W \biggl( \frac {\tau (t-[z_1^{-1}])}{\tau (t)} ,
                       \frac {\tau (t-[z_2^{-1}])}{\tau (t)}\biggr) \,  -  
   (z_1-z_2) \, \frac {\tau (t-[z_1^{-1}]) \tau (t-[z_2^{-1}])}{\tau ^2(t)}
                                                              \biggr] \, = 
$$ 
$$
{e}^{(z_1+z_2) x} \biggl[ \frac {W \bigl( \tau (t-[z_1^{-1}]),
                                  \tau (t-[z_2^{-1}]) \bigr)}{\tau ^2(t)} -   
     (z_1-z_2) \, \frac {\tau (t-[z_1^{-1}]) \tau (t-[z_2^{-1}])}{\tau ^2(t)}
                                                              \biggr] \, = 
$$
$$
\align
\frac {{e}^{(z_1+z_2) x}}{\tau ^2(t)} \biggl[ 
    (z_1-z_2) \bigl( \tau (t-[z_1^{-1}]) \tau (t-&[z_2^{-1}]) - 
                     \tau (t) \tau (t-[z_1^{-1}]-[z_2^{-1}]) \bigr)      - \\
    &(z_1-z_2) \tau (t-[z_1^{-1}]) \tau (t-[z_2^{-1}]) \biggr]           = 
\endalign
$$
$$
    (z_1-z_2) \, {e}^{x \, (z_1+z_2)} 
 \frac {\tau (t-[z_1^{-1}]-[z_2^{-1}])}{\tau (t)} \, . 
$$
From here follows (i), because we have ($\, t_1 \equiv x$ ) : 
$$
{e}^{x \, (z_1+z_2)} \exp \bigl( \sum_{k=1}^{\infty} t_{2k+1} 
                                (z_1^{2k+1}+z_2^{2k+1}) \bigr) =  
\exp \bigl( \sum_{k=0}^{\infty} t_{2k+1}(z_1^{2k+1}+z_2^{2k+1}) \bigr).  
$$
It is easy to proove (ii), (iii) and (iv) in the same way.
\qed

\bigpagebreak

\subhead
III. Proof of the main results
\endsubhead

\proclaim
{Proof of Theorem 1.2}
\endproclaim

First we will proove the identity (i).
Using the identities of Lemma 2.2 and (1.5), let us expand the following 
Wronskian :
$$
W \bigl( \tau (t+[z_1]) \tau (t-[z_1]) \, , \, 
         \tau (t+[z_2]) \tau (t-[z_2]) \bigr) \, 
$$
in two different ways.
From the first line of (1.5) we obtain : 
$$
\align
(z_2^{-1}-z_1^{-1})\tau (t) \bigl[ &\tau (t-[z_1]-[z_2]) \tau (t+[z_1]) 
                                                        \tau (t+[z_2]) - \\
   &\tau (t+[z_1]+[z_2]) \tau (t-[z_1]) \tau (t-[z_2]) \bigr] \, ,
\endalign
$$
and from the second line of (1.5) we obtain :
$$
\align
(z_2^{-1}+z_1^{-1})\tau (t) \bigl[ &\tau (t+[z_1]-[z_2]) \tau (t-[z_1]) 
                                                        \tau (t+[z_2]) - \\
   &\tau (t-[z_1]+[z_2]) \tau (t+[z_1]) \tau (t-[z_2]) \bigr] \, .
\endalign
$$
But $(z_2^{-1}-z_1^{-1}) = {(z_1-z_2)}/{z_1z_2}$ and 
 $(z_2^{-1}+z_1^{-1}) = {(z_1+z_2)}/{z_1z_2}$\, , so we have :
$$
\frac {\tau (t)}{z_1z_2} \, \biggl[\, \text{\it l.h.s. of   
(i)} \, \biggr]  \qquad = \qquad 
 \frac {\tau (t)}{z_1z_2} \, \biggl[ \, \text{\it r.h.s. of   
(i)} \,  \biggr]  \, .
$$
The proof of the first identity of Theorem 1.2 is completed.

Now we will obtain the second identity (ii) of Theorem 1.2 letting 
$\, z_2 \, $ to tend to $\, z_1 \, $   
in the first identity (i) (we will denote $\, z_1 = z_2 = z \, $) .
The {\it l.h.s.} of (ii) is clear. In order to obtain the {\it r.h.s.} 
of (ii) we mention that :
$$
\align
&\partial_{z_2} \bigl( \tau (t+[z_1]-[z_2])\left. \bigr)\right|_{z_2 = 
                                                             z_1 =z} = \\
\partial_{z_2} \biggl[ \tau \bigl( (x+z_1)-z_2 \, , \, &(t_3+\frac {z_1^3}{3}) 
- \frac{z_2^3}{3} \, , \, (t_5+\frac{z_1^5}{5})-\frac{z_2^5}{5} \, , \, 
              \dots \bigr)\left.  \biggr] \right|_{z_2 = z_1 =z}  = \\
&- \sum_{k=0}^{\infty} z^{2k} \partial_{t_{2k+1}} \tau (t) .
\endalign
$$
The same way we obtain :
$$
\align
&\partial_{z_2} \bigl( \tau (t-[z_1]+[z_2])\left.\bigr)\right|_{z_2=z_1=z} = \,
 \sum_{k=0}^{\infty} z^{2k} \partial_{t_{2k+1}} \tau (t) \, , \\
&\partial_{z_2} \bigl( \tau (t+[z_2]) \left. \bigr)\right|_{z_2=z} = \,
 \sum_{k=0}^{\infty} z^{2k} \partial_{t_{2k+1}} \tau (t+[z]) \, , \\
&\partial_{z_2} \bigl( \tau (t-[z_2]) \left. \bigr)\right|_{z_2=z} = \, - \,  
 \sum_{k=0}^{\infty} z^{2k} \partial_{t_{2k+1}} \tau (t-[z]) \, .
\endalign
$$
So, from the {\it r.h.s.} of (i) we obtain :
$$
\align
\tau (t-[z]) \, \sum_{k=0}^{\infty} 2z^{2k+1} \biggl[ \tau (t) 
                              \partial_{t_{2k+1}} \, \tau (t+[z]) \,  &- 
   \tau (t+[z]) \partial_{t_{2k+1}} \, \tau (t) \biggr] \, + \\
\tau (t+[z]) \, \sum_{k=0}^{\infty} 2z^{2k+1} \biggl[ \tau (t) 
                              \partial_{t_{2k+1}} \, \tau (t-[z]) \,  &- 
   \tau (t-[z]) \partial_{t_{2k+1}} \, \tau (t) \biggr] \, , 
\endalign
$$ 
which gives the {\it r.h.s.} of (ii).
\qed
 
\proclaim
{Proof of Theorem 1.1} 
\endproclaim

On the one hand, using the expressions of the wave functions $\psi (t,z)$ 
and $\psi ^{*} (t,z)$ in terms of tau function $\tau (t)$ (in our case 
$\tau$ is an arbitrary KdV tau function) we obtain from the first line  
of (1.4) : 
$$
\bold W \, \equiv \, 
W \bigl[ \, \psi (t,z_1)\psi ^{*}(t,z_1)\, , \,
                             \psi (t,z_2)\psi ^{*}(t,z_2) \, \bigr]       = 
$$
$$
\align
W \biggl( \, \frac {\tau (t+[z_1^{-1}])\tau (t-[z_1^{-1}])} 
                                                         {\tau ^2 (t)} \, , \, 
             \frac {\tau (t+[z_2^{-1}])\tau (t-[z_2^{-1}])} 
                                                 {\tau ^2 (t)} \, &\biggr) = \\
\frac {1}{\tau ^4 (t)} \, 
                  W \biggl( \, \tau (t+[z_1^{-1}])\tau (t-[z_1^{-1}]) \, , \, 
                    \tau (t+[z_2^{-1}])\tau (t-[z_2^{-1}]) \, &\biggr) \, .
\endalign 
$$
From the proof of the identity (i) of the Theorem 1.2 we know that 
this equals either to :
$$
\align
\frac {z_2-z_1}{\tau ^3 (t)} \, \biggl[ 
           &\tau (t-[z_1^{-1}]-[z_2^{-1}]) 
                             \tau (t+[z_1^{-1}]) \tau (t+[z_2^{-1}]) - \\
           &\tau (t+[z_1^{-1}]+[z_2^{-1}]) 
                     \tau (t-[z_1^{-1}]) \tau (t-[z_2^{-1}]) \biggr] \, ,
\endalign
$$
 or to : 
$$
\align
\frac {z_2+z_1}{\tau ^3 (t)} \biggl[
           &\tau (t+[z_1^{-1}]-[z_2^{-1}]) 
                       \tau (t-[z_1^{-1}]) \tau (t+[z_2^{-1}]) - \\
           &\tau (t-[z_1^{-1}]+[z_2^{-1}]) 
                       \tau (t+[z_1^{-1}]) \tau (t-[z_2^{-1}]) \biggr] \, .
\endalign
$$
On the other hand, using the relations from Lemma 2.1 ,
Lemma 2.2 and Lemma 2.3 we obtain from the the second line of (1.4) :
$$
\bold W _1 \equiv -(z_1^2-z_2^2)^{-1} \, \partial_{x} \bigg[ 
W \bigl[ \psi (t,z_1) , \psi (t,z_2) \bigr] \, 
W \bigl[ \psi ^{*} (t,z_1) , \psi ^{*} (t,z_2) \bigr] \, \biggr] \, = 
$$
$$
(z_1-z_2)^2 \, (z_1^2-z_2^2)^{-1} \, \partial_{x} \biggl[ \, 
\frac {\tau (t-[z_1^{-1}]-[z_2^{-1}]) \tau (t+[z_1^{-1}]+[z_2^{-1}])}  
{\tau ^2 (t)} \, \biggr]  \, =
$$
$$
\align 
- \, \frac {z_1-z_2}{z_1+z_2}  \tau ^{-3} (t)
\biggl[ \tau (t-[z_1^{-1}]-[z_2^{-1}]) W \bigl( \tau (t+[z_1^{-1}]+[z_2^{-1}])
        , \tau (t) \bigr) &+ \\
         \tau (t+[z_1^{-1}]+[z_2^{-1}]) W \biggl(  
\tau (t-[z_1^{-1}]-[z_2^{-1}]), \tau (t) \biggr) \biggr] \, &= 
\endalign
$$
$$
\align 
 &\frac {z_2-z_1}{z_1+z_2}  \tau ^{-3} (t) \\ 
  \biggl[ &\tau (t-[z_1^{-1}]-[z_2^{-1}]) \biggl( -(z_1+z_2) \bigl( 
      \tau (t+[z_1^{-1}]+[z_2^{-1}]) \tau (t) - 
      \tau (t+[z_1^{-1}]) \tau (t+[z_2^{-1}]) \bigr) \biggr)  + \\
 &\tau (t+[z_1^{-1}]+[z_2^{-1}]) \biggl( (z_1+z_2) \bigl( 
 \tau (t-[z_1^{-1}]-[z_2^{-1}]) \tau (t) - 
\tau (t-[z_1^{-1}]) \tau (t-[z_2^{-1}]) \bigr) \biggr) \biggr] =
\endalign
$$ 
$$
\align
\frac {z_2-z_1}{\tau ^3 (t)} \, \biggl[ 
           &\tau (t-[z_1^{-1}]-[z_2^{-1}]) 
                             \tau (t+[z_1^{-1}]) \tau (t+[z_2^{-1}]) - \\ 
           &\tau (t+[z_1^{-1}]+[z_2^{-1}]) 
                     \tau (t-[z_1^{-1}]) \tau (t-[z_2^{-1}]) \biggr] \, .
\endalign
$$
and for the third line of (1.4) : 
$$
\bold W_2 \equiv -(z_1^2-z_2^2)^{-1} \, \partial_{x} \bigg[ 
W \bigl[ \psi (t,z_1) , \psi ^{*} (t,z_2) \bigr] \, 
W \bigl[ \psi ^{*} (t,z_1) , \psi (t,z_2) \bigr] \, \biggr] \, =
$$
$$ 
(z_1+z_2)^2 \, (z_1^2-z_2^2)^{-1} \, \partial_{x} \biggl[ \, 
\frac {\tau (t-[z_1^{-1}]+[z_2^{-1}]) \tau (t+[z_1^{-1}]-[z_2^{-1}])}  
{\tau ^2 (t)} \, \biggr]  \, =
$$
$$
\align 
- \, \frac {z_1+z_2}{z_1-z_2}  \tau ^{-3} (t)
\biggl[ \tau (t-[z_1^{-1}]+[z_2^{-1}]) W \bigl( \tau (t+[z_1^{-1}]-[z_2^{-1}])
        , \tau (t) \bigr) &+ \\
         \tau (t+[z_1^{-1}]-[z_2^{-1}]) W \biggl(  
          \tau (t-[z_1^{-1}]+[z_2^{-1}]), \tau (t) \biggr) \biggr] \, &= 
\endalign
$$
$$
\align
&\frac {z_2+z_1}{z_1-z_2}  \tau ^{-3} (t) \\
  \biggl[ &\tau (t-[z_1^{-1}]+[z_2^{-1}]) \biggl( -(z_1-z_2) \bigl( 
      \tau (t+[z_1^{-1}]-[z_2^{-1}]) \tau (t) - 
      \tau (t+[z_1^{-1}]) \tau (t-[z_2^{-1}]) \bigr) \biggr)  + \\
 &\tau (t+[z_1^{-1}]-[z_2^{-1}]) \biggl( (z_1-z_2) \bigl( 
 \tau (t-[z_1^{-1}]+[z_2^{-1}]) \tau (t) - 
\tau (t-[z_1^{-1}]) \tau (t+[z_2^{-1}]) \bigr) \biggr) \biggr] =
\endalign
$$
$$
\align  
\frac {z_2+z_1}{\tau ^3 (t)} \, \biggl[ 
           &\tau (t+[z_1^{-1}]-[z_2^{-1}]) 
                             \tau (t-[z_1^{-1}]) \tau (t+[z_2^{-1}]) - \\ 
           &\tau (t-[z_1^{-1}]+[z_2^{-1}]) 
                     \tau (t+[z_1^{-1}]) \tau (t-[z_2^{-1}]) \biggr] \, .
\endalign
$$

This way we obtain that $\bold W$ equals either to $\bold W_1$ or to 
$\bold W_2$ , i.e. the Faddeev - Tahtajan identity is fulfiled.
\qed

\bigpagebreak

\subhead
IV. Conclusion remarks and examples
\endsubhead
  
Firstly we illustrate the identities from Theorem 1.2 by examples with 
polynomial KdV tau functions. The author thanks F.A. Gr\"unbaum 
for sugestions to include these examples in the body of the paper.

{\bf Example 4.1}
The first nontrivial polynomial KdV tau function is $\, \tau_1 (t) := t_1 $. 
In this case the examination of the identities (i) and (ii) of Theorem 1.
is easy to do directly and the result is: the both sides of (i) are 
equal to $2 z_1 z_2^3 - 2 z_1^3 z_2$, the both sides of (ii) are equal to
$4 z^3$. 

{\bf Example 4.2} 
The next polynomial KdV tau function is of degree 3: 
$\, \tau_3 (t) := t_1^3 - 3 t_3 $ and as 
is clear from the results, the examination of the identities (i) and (ii) of
Theorem 1. in this case is difficult to do directly. We used the system 
Maple V Release 4 at RIMS, Kyoto University. So, the both sides of (i)
are equal to:
$$
6 ( z_1 z_2^3 - z_1^3 z_2 ) t_1^6 + 36 ( z_1^5 z_2 - z_1 z_2^5 ) t_1^4 + 
126 ( z_1 z_2^3 - z_1^3 z_2) t_1^3 t_3  + 
$$
$$
54 ( z_1^3 z_2^5 - z_1^5 z_2^3 ) t_1^2 + 
54 ( z_1^5 z_2 - z_1 z_2^5 ) t_1 t_3 + 54 ( z_1 z_2^3 - z_1^3 z_2 ) t_3^2 \, ,
$$
and the both sides of (ii) are equal to:
$$
12 z^3 t_1^6 - 144 z^5 t_1^4 + 252 z^3 t_1^3 t_3 + 108 z^7 t_1^2 -  
216 z^5 t_1 t_3 + 108 z^3 t_3^2 \, .
$$

There were some problems with fixing the correct form of the KdV 
tau function $\, \tau_3 (t)$ - polynomial of the form $\, t_1^3 - a t_3$. 
The function $\, t_1^3 - a t_3$ satisfies the Fay identity (1.1) iff $a = 3$ . 

{\bf Remark 4.3} 
Applying the identities (1.5) to the Wronskian:
$$
W \bigl( \tau (t+[z_1]) \tau (t-[z_1]) \tau (t+[z_3]) \tau (t-[z_3]) \, , \,   
         \tau (t+[z_2]) \tau (t-[z_2]) \tau (t+[z_4]) \tau (t-[z_4]) 
 \bigr), 
$$ 
( $z_1, z_2, z_3, z_4 \in \Bbb C$ ) we can obtain 8 different 
(equivalent)
expessions where we have Wronskians of two tau functions only 
(i.e. without any Wronskian of products of tau functions).
The expressions are separated in two groups and
applying Lemma 2.2 we could see that the 
resulting identities among the expressions in each group 
are easyly obtained using the result of Theorem 1.2 (i).
The equality of the given bellow expressions (from the two groups) 
is a non-trivial seventh-order 
(specific for KdV tau functions only) identity:
$$
\align
&(z_4^{-1} -z_3^{-1})  \tau (t+[z_1]) \tau (t-[z_1]) 
                       \tau (t+[z_2]) \tau (t-[z_2]) \\
\bigl[ \tau (t+[z_3]) &\tau (t+[z_4])\tau (t-[z_3]-[z_4]) - 
       \tau (t-[z_3])  \tau (t-[z_4])\tau (t+[z_3]+[z_4]) \bigl] \, + \, \\
&(z_2^{-1} -z_1^{-1})  \tau (t+[z_3]) \tau (t-[z_3]) 
                       \tau (t+[z_4]) \tau (t-[z_4]) \\
\bigl[ \tau (t+[z_1]) &\tau (t+[z_2])\tau (t-[z_1]-[z_2]) - 
       \tau (t-[z_1])  \tau (t-[z_2])\tau (t+[z_1]+[z_2]) \bigl] \, = \, \\
&(z_2^{-1} -z_3^{-1})  \tau (t+[z_1]) \tau (t-[z_1]) 
                       \tau (t+[z_4]) \tau (t-[z_4]) \\
\bigl[ \tau (t+[z_2]) &\tau (t+[z_3])\tau (t-[z_2]-[z_3]) - 
       \tau (t-[z_2])  \tau (t-[z_3])\tau (t+[z_2]+[z_3]) \bigl] \, + \, \\
&(z_4^{-1} -z_1^{-1})  \tau (t+[z_2]) \tau (t-[z_2]) 
                       \tau (t+[z_3]) \tau (t-[z_3]) \\
\bigl[ \tau (t+[z_1]) &\tau (t+[z_4])\tau (t-[z_1]-[z_4]) - 
       \tau (t-[z_1])  \tau (t-[z_4])\tau (t+[z_1]+[z_4]) \bigl] \, .
\endalign
$$
It is clear that this way we can obtain generalized identities 
of order $2^n - 1$ for any $\, n = 4, 5, \dots $. 
The identities from Theorem 1.2 and Remark 4.3 correspond to the cases 
$n=2$ and $n=3$ respectively.  

{\bf Example 4.4} For the first polynomial tau function 
$\, \tau_1 (t) = t_1 $ the both sides of the identity from Remark 4.3 
are equal to:
$$
2(-z_1^2 +z_2^2 -z_3^2 +z_4^2)t_1^4 +4(z_1^2 z_3^2 -z_2^2 z_4^2)t_1^2 +
2(-z_1^2 z_2^2 z_3^2 +z_1^2 z_2^2 z_4^2 
  - z_1^2 z_3^3 z_4^2 +z_2^2 z_3^2 z_4^2) .
$$
For the next polynomial KdV tau function 
$\, \tau_3 (t) = t_1^3 - 3 t_3 $ 
the both sides of this identity have too many terms (more than 250). 

{\bf Remark 4.5} 
As we mentioned in the Introduction, there are some problems to translate  
the identity (i) from Theorem 1.2 to the cases when KdV tau function is 
expressed by trigonometric functions 
or elliptic theta functions. The problems 
come roughly speaking from the fact that in the original Fay identity 
(i.e. for theta functions related to Jacobians) is used the ``Prime Form'' 
(e.g. in the $g=1$ case: $\theta_{11} (z_0-z_1)
(\theta_{11}^{'} (0))^{-1}$), 
but in the Fay identity (1.1) for KP tau functions is used the difference
$(z_0-z_1)$ instead. Our next task is to fix these problems and
to find ``geometric'' explanation of the identities from the present paper.
It will be done in some next article.

{\bf Remark 4.6}
The ``elliptic version'' of the identity (ii) from Theorem 1.2 is the
following relation:
$$
\align
\theta_{11}^{'} (0) &\biggl[ \theta_{11}(x+2z) \theta_{11}^2 (x-z)- 
\theta_{11}(x-2z)\theta_{11}^2 (x+z) \biggr] =  \\ 
\theta_{11}(2z)
&\biggl[\theta_{11}(x-z) W( \theta_{11}(x), \theta_{11}(x+z) )+ 
       \theta_{11}(x+z) W( \theta_{11}(x), \theta_{11}(x-z) ) \biggr] 
\endalign
$$
It is easily obtained from the elliptic version of the identity (i) from
Theorem 1.2 (cf. the Introduction) letting $z_1 \to z_2$ and denoting 
$z_1 = z_2 \equiv z$ .

\bigpagebreak

{\bf Acknowledgements.} Firstly I would like to express my sincere thanks 
to my host researcher Prof. T. Miwa and to Prof. M. Jimbo for making it
possible to be at RIMS. 
I wish to acknowledge the useful discussions with 
F.A. Gr\"unbaum,
M. Jimbo,
T. Miwa,
M. Ohmiya, 
E. Previato, 
T. Shiota,  
B. Sturmfels. 
I am also thankfull to all people at RIMS for warm hospitality.

The work was done under JSPS Postdoctoral Fellowship for Foreign 
  Researchers.    

\bigpagebreak

\subhead
References
\endsubhead

$^{1}$ D. Mumford, ``{\it Tata lectures on theta}'' I, II, 
            Progress in Mathematics 28, 43, Birkh\"auser Boston Inc.,
            Boston, Mass., 1983, 1984.

$^{2}$ M. Adler and P. van Moerbeke, ``{\it A matrix integral solution to  
           two-dimensi- onal $W\sb p$-gravity}'', Comm. Math. Phys. 147,
           no. 1, 25--56 (1992).

$^{3}$ T. Shiota, ``{\it Characterization of Jacobian varieties in 
            terms of soliton equations}'', Invent. Math. 83, 
            no. 2, 333--382 (1986).  

$^{4}$ J. Fay, ``{\it Theta functions on Riemann surfaces}'',
            Lecture Notes in Mathematics, vol. 352, 
            Springer-Verlag, Berlin-New York, 1973.

$^{5}$ L.D. Faddeev and L.A. Takhtajan, ``{\it Hamiltonian methods in 
             the theory of solitons}'', Springer Series in Soviet Mathematics,
             Springer-Verlag, Berlin-New York, 1987. 

$^{6}$ Y.P. Mishev, ``{\it Crum-Kre\u\i n transforms and
            $\Lambda$-operators for radial Schr\"odin- ger equations}'',
            Inverse Problems 7, no. 3, 379--398 (1991).

$^{7}$ G. Wilson, ``{\it Habillage et fonctions $\tau$}'',
           C. R. Acad. Sci. Paris Ser. I Math. 299, no. 13, 
           587--590 (1984). 

$^{8}$ M. Adler and P. van Moerbeke, ``{\it Birkhoff strata, 
           B\"acklund transformations, and regularization of isospectral 
           operators}'', Adv.Math. 108, no. 1, 140-204 (1994). 

$^{9}$ G. Felder and A. Veselov, ``{\it Shift operators for the quantum 
            Calogero-Suther- land problems via Knizhnik-Zamolodchikov 
            equation}'', Commun. Math. Phys. 160, no. 2, 259-273 (1994).

$^{10}$ R. Carroll, ``{\it Some kernels on a Riemann surface}'', 
            Appl. Anal. 65, no. 3-4, 333-352 (1997). 

$^{11}$ R. Carroll and Jen Hsu Chang,    
               ``{\it The Whitham equations revisited}'', 
               Appl. Anal. 64, no. 3-4, 343-378 (1997).

$^{12}$ Y.P. Mishev, ``{\it On one cubic identity for KdV tau functions}'',
            Matematica Balkanika (1995).

$^{13}$ Y.P. Mishev, ``{\it On the Faddeev-Takhtajan identity}'',
            Matematica Balkanika (1995).

\enddocument